\documentclass[10pt]{amsart}
\usepackage[margin=1in]{geometry}
\usepackage{amsmath}
\usepackage{amsfonts,amssymb,graphicx}
\usepackage{amsthm}
\usepackage{tikz}
\usetikzlibrary{matrix,arrows,arrows.meta}
\newcommand{\ve}{\varepsilon}

\newcommand{\W}{\mathbb{W}}
\newcommand{\bs}{\mathbb{S}}
\newcommand{\F}[0]{\mathcal F}
\newtheorem{defn}{Definition}[section]
\newtheorem{prop}[defn]{Proposition}
\newtheorem{cor}[defn]{Corollary}
\newtheorem{theorem}[defn]{Theorem}
\newtheorem{lemma}[defn]{Lemma}

\graphicspath{{.}}

\usepackage{verbatim}

\DeclareMathOperator{\aut}{Aut}

\usepackage[style=alphabetic,
            url=false,
            isbn=false,
            doi=false]{biblatex}
\title{Linearity of Generalized Cactus Groups}
\author{Runze Yu}
\address{Department of Mathematics, University of California, Los Angeles}
\email{yurunze2023@ucla.edu}

\begin{document}
\maketitle

\begin{abstract}
	Cactus groups are traditionally defined based on symmetric groups, and pure cactus groups are particular subgroups of cactus groups. Mostovoy \cite{Mos} showed that pure cactus groups embed into right-angled Coxeter groups. We generalize this result to cactus groups associated with arbitrary finite Coxeter groups and we investigate some representations of generalized cactus groups and deduce the linearity of generalized cactus groups.
\end{abstract}

\section{Introduction}
For an integer $n>0,$ the \textit{cactus group} $J_n$ is the group generated by $s_{p,q}$ for $1\le p<q\le n$ with relations
\begin{align*}
&s_{p,q}^2=1, & \\
&s_{p,q}s_{m,r}=s_{m,r}s_{p,q}, & \text{if } [p,q]\cap[m,r]=\varnothing\\
&s_{p,q}s_{m,r}=s_{p+q-r,p+q-m}s_{p,q}, & \text{if } [p,q]\supset[m,r]
\end{align*}

There is a homomorphism $g:J_n\to S_n$ given by $s_{p,q}\mapsto \sigma_{p,q},$ where $\sigma_{p,q}$ reverses the order of elements $p,\cdots,q$ and leaves the rest unchanged. The kernel of $g$ is the \textit{pure cactus group} $\Gamma_n.$  

The subject of cactus groups has generated much interest. Mostovoy showed that the pure cactus group $\Gamma_n$ embeds into the diagram group $D_n,$ a right-angled Coxeter group, and $\Gamma_n$ is residually nilpotent, see \cite{Mos}. Davis, Januszkiewicz, and Scott studied the connection between cactus groups and the fundamental groups of blow-ups in \cite{DaJaSc}. Henriques and Kamnitzer introduced actions of cactus groups on tensor products of crystals and coboundary categories similar to the action of braid groups on braided categories, see \cite{HenKam}. Losev studied the case of cactus group defined from the Weyl group of a Lie algebra and the action of the cactus group on the Weyl group, see \cite{Losev}, and Bonnafé extended Losev's construction to any Coxeter group instead of a Weyl group, see \cite{Bon}. 

In a Coxeter system $(W,S),$ a nonempty subset $I\subset S$ is \emph{connected} if the subgroup $W_I$ of $W$ generated by $I$ is finite, and there is no way to write $I=J\cup K,$ where $J$ and $K$ are disjoint nonempty sets and $W_I=W_J\times W_K.$ Set $\mathcal F(S)$ to be the family of all connected subsets of $S.$ 
\begin{defn}
\normalfont
The \emph{generalized cactus group} $C_W$ is generated by $\{\gamma_{I}:I\in\mathcal F(S)\}$ subject to relations:
\begin{align*}
&\gamma_I^2=1, \\
&\gamma_I\gamma_J=\gamma_J\gamma_{w_J(I)},&\mathrm{if}\,\, I\subset J\,\,\mathrm{or} \,\,W_{I\cup J}=W_I\times W_J.
\end{align*}
\end{defn}
Here $w_J$ is the longest element in $W_J.$ There is a surjective homomorphism $g_W:C_W\to W$ given by $\gamma_I\mapsto w_I.$ The kernel of $g_W$ is the \emph{generalized pure cactus group}. Our main result is that the generalized pure cactus group can be embedded into a right-angled Coxeter group. Given a Coxeter group $(W,S)$ and $\mathcal F(S)$ the family of connected subsets of $S,$ we introduce the right-angled Coxeter group $\W$ generated by $\bs=\{w(W_I):w\in W,I\in\mathcal F(S)\}$ associated with a Coxeter matrix
$$m(W',W'')=\begin{cases}
1,& \text{if } W'=W'' \\
2,& \text{if } W'\subset W'' \text{ or } W''\subset W'	\\
2,& \text{if } W'\cap W''=\{1\} \text{ and } W''\subset C_W(W') \\
\infty, &\text{otherwise}.
\end{cases}$$
In this setting, the group $W$ acts on $(\W,\bs)$ by $w\mapsto g_w,$ where $g_w(\tau_{W'})=\tau_{wW'w^{-1}}.$ We prove the generalized cactus group $C_W$ defined on $(W,S)$ embeds into $\W\rtimes\aut(\W,\bs),$ and thereby show that
\begin{theorem}
\normalfont
	The generalized pure cactus group defined by the Coxeter system $(W,S)$ embeds into the right-angled Coxeter group $(\W,\bs).$
	\label{main}
\end{theorem}
Finally, we construct two representations of the generalized cactus group. The generalized cactus group is a special case of the group $A$ of all lifts of the $W$-action on the blow-up $\Sigma_{\#}$ of a reflection tiling $\Sigma$ to its universal cover \cite{DaJaSc}. We adapt the representation of $A$ in \cite{DaJaSc} to construct a representation of generalized cactus groups. 

We construct another representation of generalized cactus groups from the geometrical representation of Coxeter groups in light of Theorem \ref{main}. Furthermore, this representation restricts to a faithful representation on the generalized cactus group. This shows
\begin{theorem}
	\normalfont
	Generalized cactus groups are linear groups.
\end{theorem}

\section{Generalized Cactus Groups}
\subsection{Definition}
The notion of cactus group can be generalized by replacing symmetric groups with general Coxeter groups. Geometrically, it is a special case of the group $A$ introduced by David, Januszkiewicz and Scott as the group of all lifts of the $W$-action on the blow-up $\Sigma_{\#}$ of a reflection tiling $\Sigma$ to its universal cover in \cite{DaJaSc}. 

Let $(W,S)$ be a finitely generated Coxeter group and $(m_{st})_{s,t\in S}$ be its associated Coxeter matrix. A \textit{spherical subgroup} of $W$ is a finite subgroup $W_I$ generated by some $I\subset S.$ We denote its longest element by $w_I$ following \cite{Hum}. A nonempty subset $I\subset S$ is \textit{connected} if $W_I$ is finite and there is no way to write $I=J\cup K$ such that $J$ and $K$ are nonempty and $W_{I}=W_J\times W_K.$ Let $\mathcal F(S)\subset\mathcal P(S)$ denote the family of all connected subsets of $S.$\
\begin{defn}
\normalfont
The generalized cactus group $C_W$ is generated by $\{\gamma_{I}:I\in\mathcal F(S)\}$ subject to relations:
\begin{align*}
&\gamma_I^2=1 \\
&\gamma_I\gamma_J=\gamma_J\gamma_{w_J(I)}&\mathrm{if}\,\, I\subset J\,\,\mathrm{or} \,\,W_{I\cup J}=W_I\times W_J.
\end{align*}

\end{defn}
Since $w_J$ is the longest element of $W_J$, it follows that $w_J(I)$ is still a subset of $S$ for $I\subset J.$ In particular if $W_{I\cup J}=W_I\times W_J,$ then $w_J(I)=I$ and $\gamma_I$ commutes with $\gamma_J.$

Similar to the case of cactus group, there is a homomorphism $g_W:C_W\to W$ given by $\gamma_I\mapsto w_I.$ Indeed, $w_I^2=1$ since it is the longest element in $W_I$. If $W_I\times W_J=W_{I\cup J},$ then $w_I$ commutes with $w_J$. When $I\subset J,$ there is $w_J(I)\subset J$ and the longest element in $W_{w_J(I)}=w_J(W_I)$ is $w_{w_J(I)}=w_Jw_Iw_J.$ Call the kernel of $g_W$ the \textit{generalized pure cactus group} $PC_W.$

\begin{prop}
\normalfont
When $W$ has type $A_{n-1},$ there is an isomorphism $C_W\simeq J_n.$
\end{prop}
\begin{proof}
When $W$ has type $A_{n-1},$ it is isomorphic to the Coxeter system $(S_n,S)$ where $S=\{s_i:1\le i\le n-1\}$ associated by an $(n-1)\times(n-1)$ Coxeter matrix $(m_{ij})$ given by
$$m_{ij}=\begin{cases}
1,&i=j, \\
2,&|i-j|\ge 2,\\
3,&|i-j|=1.
\end{cases}$$
Identify subsets of $\{1,\cdots,n-1\}$ with subsets of $S$ via $I\mapsto \{s_i:i\in I\}$ and let $C_n:=C_W$ be the generalized cactus group on $(W,S).$ If $J$ and $K$ are subsets of $\{1,2,\cdots,n-1\}$ where $x-y\ge2$ for all $x\in J$ and $y\in K,$ then $W_J\cap W_K=\varnothing$ and every generator of $W_J$ commutes with every generator in $W_K.$ Hence $J\cup K\notin \mathcal F(S)$ and it follows that $\mathcal F(S)$ consists of connected intervals $[p,q]$ where $1\le p\le q\le n-1.$ 

Let $J_n$ be the cactus group generated by $s_{p,q}$ for $1\le p<q\le n.$ We claim that the map $\varphi:J_n\to C_n$ given by $s_{p,q}\mapsto \gamma_{[p,q-1]}$ is a group isomorphism.

We show first that $\varphi$ is a homomorphism. Firstly there is $\gamma^2_{[p,q-1]}=1.$ If $[p,q]\cap[m,r]=\varnothing,$ then by the above discussion,
$$\gamma_{[p,q-1]}\gamma_{[m,r-1]}=\gamma_{[m,r-1]}\gamma_{[p,q-1]}.$$
If $[p,q]\supset[m,r],$ then the longest element generated by $s_p,\cdots,s_{q-1}$ is $\sigma_{p,q}.$ Then
$$\gamma_{[m,r-1]}\gamma_{[p,q-1]}=\gamma_{[p,q-1]}\gamma_{\sigma_{p,q}([m,r-1])}=\gamma_{[p,q-1]}\gamma_{[p+q-r,p+q-m-1]}.$$
Conjugating by $\gamma_{[p,q-1]}$ gives the corresponding relations in $J_n$ and therefore $\varphi$ is indeed a group homomorphism.

Define $\psi:C_n\to J_n$ by $\gamma_{[p,q-1]}\mapsto s_{p,q}.$ We have $s_{p,q}^2=1.$ If $W_{[p,q-1]}\times W_{[m,r-1]}=W_{[p,q-1]\cup [m,r-1]},$ then every element in $[p,q-1]$ commutes with every element in $[m,r-1].$ Since they are connected intervals, this implies $[p,q]\cap [m,r]=\varnothing$ and then $s_{p,q}s_{m,r}=s_{m,r}s_{p,q}.$ If $[p,q-1]\supset [m,r-1]$, then the relation in $C_n$ 
corresponds to the relation $s_{p,q}s_{m,r}=s_{p+q-r,p+q-m}s_{p,q}$ in $J_n,$ which is valid. It follows that $\psi$ is also a group homomorphism. Hence $C_n$ is isomorphic to $J_n$ because $\varphi\circ\psi=\mathrm{id}_{C_n}$ and $\psi\circ\varphi=\mathrm{id}_{J_n}.$
\end{proof}

\subsection{Embedment of the Generalized Pure Cactus Group}
Our main theorem in this section is that the generalized pure cactus group can be embedded into a right-angled Coxeter group, generalizing the approach taken by Mostovoy in \cite[][Proposition 2]{Mos}:
\begin{theorem}
\normalfont
Suppose $C_W$ is the generalized cactus group defined on the Coxeter system $(W,S).$ Then the generalized pure cactus group $PC_W$ embeds into a right-angled Coxeter group.
\label{embed}
\end{theorem}

We set up a technical lemma before giving the proof. Let $(W,S)$ be a right-angled Coxeter group and $G$ a subgroup of $\aut(W,S).$ Let $I$ be a subset of $S$ and $\{g_i\}_{i\in I}$ a family of involutions in $G$ satisfying
\begin{enumerate}
\item $g_i(s)\in S$ for all $i\in I$ and $s\in S.$
\item Given $i,j\in I$ with $ij=ji,$ either
$$g_i(j)\in I, g_j(i)=i,\text{ and } g_{g_i(j)}=g_ig_jg_i,$$
or
$$g_j(i)\in I, g_i(j)=j,\text{ and } g_{g_j(i)}=g_jg_ig_j.$$
\item $g_i(i)=i$ for all $i\in I.$
\end{enumerate}
Let $H$ be the subgroup of $W\rtimes G$ generated by $\{ig_i\}_{i\in I}.$ Let $L$ be the group generated by $l_i$ for all $i\in I$ subject to relations
\begin{enumerate}
\item $l_il_j=l_{g_i(j)}l_i$ when $ij=ji, g_j(i)=i,g_i(j)\in I,$ and $g_{g_i(j)}=g_ig_jg_i,$
\item $l_i^2=1$ for all $i\in I.$
\end{enumerate}
\begin{lemma}
\normalfont
The map $\varphi:L\to H$ defined by $l_i\mapsto ig_i$ is a group isomorphism.
\label{L-group}
\end{lemma}
\begin{proof}
We show first $\varphi$ is a group homomorphism. Because $g_i$ are involutions, there is $(ig_i)^2=i^2g_i^2=1.$ Now suppose $ij=ji,g_i(j)\in I, g_j(i)=i,$ and $g_{g_i(j)}=g_ig_jg_i.$ Then
$$(g_i(j)g_{g_i(j)})(ig_i)=(g_i(j)(g_ig_jg_i)(i))(g_ig_jg_ig_i)=g_i(ji)(g_ig_j)=g_i(ij)(g_ig_j)=(ig_i(j))(g_ig_j)=(ig_i)(jg_j).$$
Hence $\varphi$ is indeed a group homomorphism. It is also surjective as $\{ig_i\}_{i\in I}$ generates $H.$ To show injectivity, let $i_1,i_2,\cdots,i_m\in I$ and let $i_r'=(g_{i_{1}}g_{i_{2}}\cdots g_{i_{r-1}})(i_r)\in S.$ Note that $i_r'$ only depends on $i_k$ for $1\le k\le r.$ We shall give a lemma first:
\begin{lemma}
	\normalfont
	Let $l_{i_1}\cdots {l_{i_m}}$ be a word in $L.$ 
	\begin{enumerate}
		\item If there is an $r$ such that $i_r'=i_{r+1}',$ then there is another word $l_{j_1}l_{j_2}\cdots l_{j_{m-2}}=l_{i_1}\cdots {l_{i_m}}$ such that $j_k'=i_k'$ for $1\le k\le r-1$ and $j_k'=i_{k+2}'$ for $r\le k\le m-2.$
		\item If $i_r'$ commutes with $i_{r+1}',$ then there is another word $l_{j_1}l_{j_2}\cdots l_{j_{m-2}}=l_{i_1}\cdots {l_{i_m}}$ such that $j_k'=i_k'$ for $k\ne r,r+1, j_r'=i_{r+1}'$ and $j_{r+1}'=i_r'.$
	\end{enumerate}
	\label{L-word}
\end{lemma}
\begin{proof}
Assume $i_r'=i_{r+1}'.$ Then $i_r=g_{i_r}(i_{r+1})$ and then $i_r=i_{r+1}$ since the $g$'s are bijective involutions and $i_r$ is fixed by $g_{i_{r+1}}.$ Then take $j_k=i_k$ for $1\le k\le r-1$ and $j_k=i_{k+2}$ for $r\le k\le m-2.$ Then  $$l_{j_1}l_{j_2}\cdots l_{j_{m-2}}=l_{i_1}\cdots l_{i_{r-1}}l_{i_{r+2}}\cdots l_{i_m}.$$
Moreover, $j_k'=i_k'$ for $1\le k\le r-1$ and $j_k'=(g_{i_1}g_{i_2}\cdots g_{i_{r-1}}g_{i_{r+2}}\cdots g_{i_{k-1}})(i_{k+2})=i'_{k+2}$ since $g_{i_r}=g_{i_{r+1}}$ and they are involutions.

Now assume $i_r'$ and $i_{r+1}'$ commute and then $i_rg_{i_r}(i_{r+1})=g_{i_r}(i_{r+1})i_r.$ Hence $i_r$ commutes with $i_{r+1}$ after applying $g_{i_r}.$ Then either $g_{i_r}(i_{r+1})\in I$ or $g_{i_{r+1}}(i_r)\in I.$ If the former is true, then $l_{i_r}l_{i_{r+1}}=l_{g_{i_r}(i_{r+1})}l_{i_r}$ and take $$l_{j_1}l_{j_2}\cdots l_{j_m}=l_{i_1}\cdots l_{i_{r-1}}l_{g_{i_r}(i_{r+1})}l_{i_r}\cdots {l_{i_m}}.$$ Then there are $j_k'=i_k'$ for $1\le k\le r-1,$ and $j_k'=i_k'$ for $k\ge r+2$ since $g_{g_{i_r}(i_{r+1})}g_{i_r}=g_{i_r}g_{i_{r+1}}.$ Finally, there is $j_r'=i_{r+1}'$ and $j_{r+1}'=j_r'$ since $g_{g_{i_r}(i_{r+1})}g_{i_r}=g_{i_r}g_{i_{r+1}}$ and $i_r$ is fixed under $g_{i_{r+1}}.$ Similarly, if the latter is true, take $$l_{j_1}l_{j_2}\cdots l_{j_m}=l_{i_1}\cdots l_{i_{r-1}}l_{i_{r+1}}l_{g_{i_{r+1}}(i_r)}{l_{i_m}}.$$
This proves Lemma \ref{L-word}.
\end{proof}
\begin{comment}

In addition, we cite the following lemma from \cite[][Theorem 3.4.2(ii)]{Davis}:
\begin{lemma}
\normalfont
Let $(W,S)$ be a right-angled Coxeter system. Two reduced expressions $w$ and $w'$ represent the same element of $W$ if and only if one can be transformed into the other by a sequence of operations that replace a subword of the form $(s, t)$ by the word $(t, s)$ where $st=ts.$
\label{davis}
\end{lemma}\end{comment}
Now suppose $\varphi(l_{i_1}l_{i_2}\cdots l_{i_m})=e.$ Then $i_1'i_2'\cdots i_n'=e.$ Since $i_1'i_2'\cdots i_n'$ is a word in the right-angled Coxeter system $(W,S),$ there is a sequence of tuples of elements of $S$
$$K^1=(K^1_1,K^1_{2},\cdots,K^1_{n_d}),\cdots,K^d=(K^d_{1},K^d_{2},\cdots,K^d_{n_d})$$ 
 such that 
 $$K^1=(K^1_1,K^1_{2},\cdots,K^1_{n_d})=(i_1',i_2',\cdots,i_m')$$
 and $K^d$ is the empty tuple representing $e,$ such that for each $s,$ either there exists an $r$ such that $K^s_{r}=K^s_{r+1}$ and 
\begin{equation}K^{s+1}=(K^s_{1},\cdots,K^s_{r-1},K^s_{r+2},\cdots,K^s_{n_s}) \label{relA}
\end{equation}
or there exists a $t$ such that $K^s_{t}$ and $K^s_{t+1}$ commute and 
\begin{equation}K^{s+1}=(K^s_{1},\cdots,K^s_{t+1},K^s_{t},\cdots,K^s_{n_s}). \label{relB}
\end{equation}
Analyze $s=1$ as an example. If (\ref{relA}) holds for $s=1,$ then $i_r'=i_{r+1}'$ for some $r$ and by part 1 of Lemma \ref{L-word}, $K^2$ corresponds to some $l_{j_1}l_{j_2}\cdots l_{j_{m-2}}=l_{i_1}\cdots {l_{i_m}}.$ If (\ref{relB}) holds for $s=1,$ then $i_r'$ and $i_{r+1}'$ commute for some $r,$ and $K^2$ corresponds to some $l_{j_1}l_{j_2}\cdots l_{j_{m}}=l_{i_1}\cdots {l_{i_m}}$ by part 2 of Lemma \ref{L-word}. By using the same argument inductively, we can conclude that $l_{i_1}\cdots {l_{i_m}}=e.$
Hence $\varphi$ is an isomorphism. This proves Lemma \ref{L-group}.
\end{proof}

Now we give the proof of Theorem \ref{embed}:
\begin{proof}
Define $\mathbb S:=\{w(W_I):I\in \mathcal \F(S) \text{ and } w\in W\}.$
Let $\mathbb M$ be the Coxeter matrix whose entries are given by $m(W',W'')\in\mathbb N\cup\{\infty\}$ for each pair $(W',W'')\in\bs\times\bs$ where
$$
m(W',W'')=\begin{cases}
1,& \text{if } W'=W'' \\
2,& \text{if } W'\subset W'' \text{ or } W''\subset W'	\\
2,& \text{if } W'\cap W''=\{1\} \text{ and } W''\subset C_W(W') \\
\infty, &\text{otherwise} 
\end{cases}
$$
Let $(\mathbb{W},\bs)$ be the right-angled Coxeter group on generators $\tau_{W'}$ associated to the Coxeter matrix $\mathbb M.$ For each $w\in W,$ there is a group automorphism $g_w\in\aut (\W,\bs)$ generated by $\tau_{W'}\mapsto \tau_{wW'w^{-1}},$ since conjugation does not affect inclusion relations between subsets of $W$ and commutativity relations between elements of $W.$ Given $I\in\mathcal F(S),$ put $g_I=g_{W_I}=g_{w_I}.$ Note that $g_{I}$ is an involution as $w_I$ is.

\begin{lemma}
\normalfont
	Let $I,J\in \mathcal F(S).$ We have
	\begin{enumerate}
		\item[(1)] $g_I(\tau_{W_I})=\tau_{W_I}.$
		\item[(2)] If $I\subset J,$ then $g_I(\tau_{W_J})=\tau_{W_J},g_J(\tau_{W_I})=\tau_{W_{w_J(I)}},$ and $g_{w_J(I)}=g_Jg_Ig_J.$
		\item[(3)] If $W_{I\cup J}=W_I\times W_J,$ then $g_I(\tau_{W_J})=\tau_{W_J},g_J(\tau_{W_I})=\tau_{W_I},g_{w_J(I)}=g_Jg_Ig_J,$ and $g_{w_I(J)}=g_Ig_Jg_I.$
		\item[(4)] $m(W_I,W_J)=2, g_J(\tau_{W_I})=\tau_{W_I}$ and $g_{w_I(J)}=g_Ig_Jg_I$ if and only if $J\subset I$ or $W_I\times W_J=W_{I\cup J}$ or $I\subset J$ with $w_J(I)=I.$
	\end{enumerate}
\end{lemma}
\begin{proof}
(1) is clear. For (2), assume $I\subset J.$ Then $w_I\in W_I\subset W_J$ and then $w_I(W_J)= W_J$ since conjugation is bijective. Hence $g_I(\tau_{W_J})=\tau_{W_J}.$ Furthermore, because $W_I$ is generated by $I,$ a generating set of $w_J(W_I)$ is $w_J(I)\subset w_J(J)=J$ since $w_J$ is the longest element. Then $g_J(\tau_{W_I})=\tau_{W_{w_J(I)}}.$ Finally, $g_{w_J(I)}$ is conjugation by the longest element in $w_J(W_I),$ which is $w_J(w_I)=w_Jw_Iw_J.$  It follows immediately that $g_{w_J(I)}=g_Jg_Ig_J.$

(3) can be shown by two similar arguments in which $w_J(W_I)=W_I$ and $w_I(W_J)=W_J$ respectively.

Sufficiency in (4) follows from (2) and (3). Now assume $m(W_I,W_J)=2, g_J(\tau_{W_I})=\tau_{W_I}$ and $g_{w_I(J)}=g_Ig_Jg_I.$ The only nontrivial case is $I\subset J$ by definition of the Coxeter matrix $M.$ In this case, (2) implies $W_I=W_{w_J(I)}=w_J(W_I).$ In particular, $w_J(I)=w_J(W_I)\cap S=W_I\cap S=I.$ This completes the proof.
\end{proof}

We conclude the proof of theorem \ref{embed} by showing the following lemma:
\begin{lemma}
\normalfont
	There is a monomorphism $C_W\hookrightarrow \W\rtimes \aut(\W,\bs)$ defined by $\gamma_I\mapsto \tau_{W_I}\cdot g_I.$ It restricts to a monomorphism $PC_W\hookrightarrow \W.$
	\label{conclusion}
\end{lemma}
\begin{proof}
There is a monomorphism $h:L\to \W\rtimes \aut(\W,\bs),$ where $L\simeq C_W$ is the group generated by $l_I$ as in Lemma \ref{L-group}. This proves the first part.

To show the second part, observe that the following diagram commutes:
\begin{center}
\begin{tikzpicture}
	\matrix (m)
		[
			matrix of math nodes,
			row sep    = 3em,
			column sep = 4em
		]
		{
			C_W  & \W\rtimes \aut(\W,\bs) \\
			W   & \aut(\W,\bs) \\
		};
	\path
		(m-1-1) edge [->] node [left] {$g_W$} (m-2-1)
		(m-1-1) edge [->] node [above] {$h$} (m-1-2)
		(m-2-1) edge [->] node [above] {$g$} (m-2-2)
		(m-1-2) edge [->>] node [right] {Proj} (m-2-2);
\end{tikzpicture}
\end{center}
in which
$$\mathrm{Proj}(h(\gamma_I))=\mathrm{Proj}(\tau_{W_I}\cdot g_I)=g_I,$$
$$g(g_W(\gamma_I))=g(w_I)=g_{w_I}=g_I.$$
Because $h$ is injective, the kernel of $g_W$ lies inside the kernel of the projection, which is isomorphic to $\W.$ This completes the proof.
\end{proof}
Theorem \ref{embed} now follows immediately from Lemma \ref{conclusion}.
\end{proof}
As an example, the usual cactus group $J_n$ defined on the Coxeter group $S_n$ coincides with the corresponding generalized cactus group $C_n.$ It follows from Theorem \ref{embed} that $J_n$ embeds into $D_n\rtimes S_n$ where $D_n$ is the diagram group and $C_n$ embeds into $\W\rtimes\aut(\W,\bs)$ for an appropriate Coxeter system $(\W,\bs),$ see \cite{Mos}. Furthermore, there is a bijection $f$ from $\bs$ to the set of subsets $I\subset\{1,2,\cdots, n\}$ with $|I|\ge2$ given by 
$$W_{\{s_i,s_{i+1},\cdots,s_{j}\}}\mapsto \{i,i+1,\cdots,j+1\}$$
and
$$wW_Iw^{-1}\mapsto w(f(W_I)).$$ 

It follows that there is a group isomorphism $\W\simeq D_n$ and a commutative diagram
\begin{center}
\begin{tikzpicture}
	\matrix (m)
		[
			matrix of math nodes,
			row sep    = 3em,
			column sep = 4em
		]
		{
			C_W  & \W\rtimes \aut(\W,\bs) \\
			J_n   & D_n\rtimes S_n \\
		};
	\path
		(m-2-1) edge [->] node [left] {$\sim$} (m-1-1)
		(m-1-1) edge [right hook->] node [above] {} (m-1-2)
		(m-2-1) edge [right hook->] node [above] {} (m-2-2)
		(m-2-2) edge [right hook->] node [right] {} (m-1-2);
\end{tikzpicture}
\end{center}

\section{Representations of Generalized Cactus Groups}
\subsection{Geometrical Representation of Coxeter Groups}
Let $(W,S)$ be a Coxeter system and $M$ its associated Coxeter matrix. A representation of $W\rtimes \aut(W,S)$ can be constructed as follows. Let $E=\mathbb R^S$ with standard basis $\{\ve_s\}_{s\in S}.$ Define a bilinear form $B_t$ on $E$ by
$$
B_t(\ve_{s},\ve_{v})=\begin{cases}
	-\cos\left(\pi/m_{s,v}\right) & m_{s,v}<\infty \\
	-t &m_{s,v}=\infty
\end{cases}
$$
For any $t>0,$ there is a representation $\pi$ of $W$ (in particular the Tits representation when $t=1$) generated by $s\mapsto \sigma_s$ where
$$\sigma_s(x)=x-2B_t(x,\ve_s)\ve_s.$$
To show $\pi$ is a representation, it suffices to show $\sigma_s\sigma_v$ has order $m_{s,v}$ for $s,v\in S.$ Clearly $\sigma_s$ preserves the bilinear form $B_t.$ Let $\lambda=B_t(e_s,e_v).$ Then
$$(\sigma_s\sigma_v)(e_s)=\sigma_s(e_s-2\lambda e_v)=-e_s-2\lambda (e_v-2\lambda e_s)=(4\lambda^2-1)e_s-2\lambda e_v$$
and $(\sigma_s\sigma_v)(e_v)=-e_v+2\lambda e_s.$ If $m_{s,v}=\infty,$ then $\lambda=-t.$ For all $t\ge1,$ the equation $k^2+2\lambda k+1=0$ has a positive real root $k=k_0>0.$ Set $r=4\lambda^2+2k_0\lambda-1$ and assume $r\ne0$ (this is true for all but finitely many $t$). Then
$$(\sigma_s\sigma_v)(e_s+k_0e_v)=(4\lambda^2+2\lambda k_0-1)e_s-(2\lambda+k_0)e_v=r(e_s+k_0e_v).$$
Hence using the second equation which can be rewritten into $(\sigma_s\sigma_v)(e_v)=2\lambda(e_s+k_0e_v)-(2\lambda k_0+1)e_v,$ we can obtain
$$(\sigma_s\sigma_v)^m(e_v)=a_m(e_s+k_0e_v)+(-2\lambda k_0-1)^me_v$$
where $a_1=2\lambda$ and $a_{m+1}=ra_m+2\lambda(-2\lambda k_0-1)^m.$ If $t>1,$ then $-2\lambda k_0-1=k_0^2>1$ and $(\sigma_s\sigma_v)^m$ can never be the identity. If $t=1,$ then $a_m=2m$ and $(\sigma_s\sigma_v)^m$ cannot be the identity, either. Thus $\sigma_s\sigma_v$ has order $\infty.$

When $m_{s,v}<\infty,$ we use the identification $\text{span}(e_s,e_v)\approx \mathbb C$ by identifying $e_s$ to 1 and $e_v$ to $-e^{-i\theta},$ where $\theta=\pi/m_{s,v}.$ Then $\sigma_s\sigma_v$ acts by a rotation of $2\pi/m_{s,v}$ and thus has order $m_{s,v}.$ Hence $\pi$ is a representation.

Moreover, a faithful representation $\pi'$ of $\aut(W,S)$ on $GL(E)$ is given by  $\pi'(g)(\ve_s)=\ve_{g(s)}$ for all basis vectors $\ve_s.$ It is possible to construct a representation of $W\rtimes \aut(W,S)$ from $\pi$ and $\pi'.$ 
 \begin{prop}
 \normalfont
Let $\Pi:W\rtimes \aut(W,S)\to GL(E)$ be given by $\Pi(\tau  g)=\pi(\tau)\pi'(g).$ Then $\Pi$ is a representation of $W\rtimes \aut(W,S).$
\end{prop}
 \begin{proof}
Since $g$ is an automorphism and preserves the Coxeter matrix, $\pi'(g)$ preserves the bilinear form $B_t.$ As $\pi$ and $\pi'$ are representations of $W$ and $\aut(W,S)$ respectively, we need to check $\Pi$ respects the product structure induced by the semi-direct product. It suffices to check that for any $g\in\aut(W,S)$ and $v\in S,$  
$$\pi'(g)\pi(v)=\pi(g(v))\pi'(g).$$
In fact, for any basis vector $\ve_{u},$ there is
\begin{align*}
(\pi'(g)\circ\pi(v))(\ve_u)&=\pi'(g)(\ve_u-2B_t(\ve_u,\ve_v)\ve_v) \\
&=\ve_{g(u)}-2B_t(\ve_{g(u)},\ve_{g(v)})\ve_{g(v)} \\
&=\pi(\tau_{g(v)})(\ve_{g(u)}) \\
&=(\pi(\tau_{g(v)})\circ\pi'(g))(\ve_u).
\end{align*}
Hence $\Pi$ is a representation.  
\end{proof}
When $t=1,$ the bilinear form is the canonical bilinear form $B$ induced by the Coxeter matrix $M.$ The following result is classical but we could not find a reference. The injectivity of $\pi$ is for example in \cite[][Corollary 5.4]{Hum}.
\begin{lemma}
    \normalfont 
    When $t=1,$ the representation $\Pi$ is faithful.
\end{lemma}
\begin{proof}
    Suppose $\Pi(wg)=\mathrm{id}.$ If $w\ne e,$ then there exists $s\in S$ such that $l(ws)<l(w).$ Let $s'=g^{-1} (s).$ Then \cite[][Theorem 5.4]{Hum} shows that $B(w(\ve_s),\ve_t)\le 0$ for all $t\in S.$ Then $B(wg(\ve_{s'}),\ve_t)\le 0.$ But $B(wg(\ve_{s'}),\ve_{s'})=B(\ve_{s'},\ve_{s'})=1,$ a contradiction. Hence $w=e$ and then $g=\mathrm{id}$ since $\pi'$ is faithful. This shows faithfulness of $\Pi.$
\end{proof}

Now pick a $t$ so that $B_t$ is nondegenerate. Theorem 1 shows that $C_W$ can be embedded into the semi-direct product $\W\rtimes \aut(\W,\bs).$ Define $\mathbb E$ to be the vector space on basis $\{\ve_{W'}\}_{W'\in\mathbb S}$ over $\mathbb R.$ The space $E$ embeds into $\mathbb E$ through the map $\ve_I\mapsto\ve_{W_I}.$ Through this map, the bilinear form $B_t$ on $E$ agrees with the bilinear form $\bar B_t$ on $\mathbb E$ generated by
$$
\bar{B_t}(\ve_{W'},\ve_{W''})=\begin{cases}
	-\cos\left(\pi/m(W',W'')\right),& m(W',W'')<\infty \\
	-t.&m(W',W'')=\infty
\end{cases}
$$
Let $\Pi$ be the representation of $\W\rtimes \aut(\W,\bs)$ on $\mathbb E$ constructed previously. By restricting to $C_W$ and its isomorphic image in $\W\rtimes \aut(\W,\bs),$ $\Pi$ becomes a representation of $C_W.$
\begin{cor}
\normalfont
The representation of $C_W$ on $\mathbb E$ is faithful. Therefore, generalized cactus groups are linear groups.
\end{cor}

\subsection{Another Linear Representation}
This section is motivated by sections 5.10, 6.1, and 6.2 in \cite{DaJaSc}. 
The generalized cactus group and certain representations can be constructed from a geometrical viewpoint as in \cite{DaJaSc}. Here we give a direct algebraic construction.
Let $(W,S)$ be a Coxeter system and let $C_W$ be the generalized cactus group on $(W,S).$ For $I,J\in \mathcal F(S),$ set 
$$
m(I,J)=\begin{cases}
1,& \text{if } I=J \\
2,& \text{if } I\subset J \text{ or } J\subset I	\\
2,& \text{if } W_I \times W_J=W_{I\cup J}\\
\infty, &\text{otherwise.}
\end{cases}
$$
Let $M$ be the Coxeter matrix indexed by $\mathcal F(S)\times\mathcal F(S)$ with entries given by $m(I,J).$ A representation of $C_W$ can be constructed as in \cite{DaJaSc}
. Let $E$ be the real vector space spanned by a basis $\{\ve_I\}_{I\in\mathcal{F}(S)}.$ Define
$$\mathcal{C}(I):=\{J\in \mathcal{F}(S):J\subset I\}.$$ 
Note that $I\notin\mathcal C(I).$ Let $E_I\subset E$ denote the subspace spanned by the set
$$\{\ve_J-\ve_{J'}: J\in\mathcal{C}(I), J'=w_I(J)\}.$$
For each $t\in \mathbb R,$ define a symmetric bilinear form on $E$ by
$$
B_t(\ve_I,\ve_J)=\begin{cases}
	1,& I=J \\
	0, &m(I,J)=2 \\
	-t,&\text{otherwise.}
\end{cases}
$$

$B_0$ is the standard inner product on $E,$ and $B_1$ is the canonical inner product associated to the Coxeter matrix $M.$ It follows that $\det B_0=1\ne 0$ and $\det B_t$ is a nonzero polynomial in $t.$ Hence $B_t$ is nondegenerate on $E$ for sufficiently large $t.$ Let $F_I$ be the orthogonal complement of $\mathbb{R}\ve_I\oplus E_I.$ Then we have an orthogonal decomposition
$$E=\mathbb{R}\ve_I\oplus E_I\oplus F_I.$$
Now define $\rho_I:E\to E$ by 
$$\rho_I=-\left.\mathrm{id}\right|_{\mathbb{R}\ve_I}\oplus-\left.\mathrm{id}\right|_{E_I}\oplus\left.\mathrm{id}\right|_{F_I}.$$
\begin{prop} 
\normalfont
	The map $\gamma_I\to\rho_I$ extends to a representation of $C_W.$ 
\end{prop}
\begin{proof}
	It suffices to show the relations in $C_W$ still hold with $\gamma_I$ replaced by $\rho_I:$
\begin{enumerate}
\item[(a)] $\rho_I^2=1,$
\item[(b)] $\rho_I\rho_J=\rho_J\rho_I$ if $W_{I\cup J}=W_I\times W_J,$ and
\item[(c)] $\rho_I\rho_J=\rho_J\rho_{w_J(I)}$ if $I\subset J.$
\end{enumerate}
$\rho_I^2=1$ since $\rho_I$ is an involution. To show (b) and (c), we give the following lemma:
\begin{lemma}
\normalfont
	\begin{enumerate}
		\item[(1)] If $W_{I\cup J}=W_I\times W_J,$ then $\rho_I(\ve_J)=\ve_J$ and $\left.\rho_I\right|_{E_J}=\mathrm{id}.$
		\item[(2)] If $J\in\mathcal{C}(I),$ then $\rho_I(\ve_J)=\ve_{J'}$ and $\rho_I(E_J)=E_{J'}$ where $J'=w_I(J).$
	\end{enumerate}
\end{lemma}
\begin{proof}
For (1), note first that $\ve_J$ is orthogonal to $\ve_I.$ Moreover, for any $L\in\mathcal C(I),$ we have $W_L\times W_J=W_{L\cup J}$ and the same is true with $L$ replaced by $w_I(L).$ Hence $\ve_J$ is orthogonal to $E_I.$ Then $\ve_J\in F_I$ and $\ve_J$ is fixed by $\rho_I.$  

Furthermore, suppose $K\in \mathcal{C}(J).$ Then $W_K\times W_I=W_{K\cup I}$ and $\rho_I(\ve_K)=\ve_K$ by the previous argument. Moreover, $K'=w_J(K)$ is also a subset of $J,$ and then $\rho_I(\ve_{K'})=\ve_{K'}$ by the same argument with $K$ replaced by $K'.$ This implies both $\ve_K$ and $\ve_{K'}$ are fixed by $\rho_I$ and then $\rho_I$ is identity on $E_J.$

For (2), if $J\in\mathcal{C}(I),$ then $J'=w_I(J)\in\mathcal C(I)$ and the vector $\ve_J+\ve_{J'}$ is orthogonal to $\ve_I.$ Let $L\in\mathcal C(I)$ and $L'=w_I(L).$ We have $\langle \ve_J,\ve_L\rangle=\langle \ve_{J'},\ve_{L'}\rangle$ and  $\langle \ve_J,\ve_{L'}\rangle=\langle \ve_{J'},\ve_{L}\rangle.$ It follows that
$$\langle \ve_J+\ve_{J'},\ve_L-\ve_{L'}\rangle=0.$$
Hence $\ve_J+\ve_{J'}$ is orthogonal to $E_I.$ Then $\ve_J+\ve_{J'}\in F_I$ and
$$\rho_I(\ve_J)=\dfrac12\left(\rho_I(\ve_J+\ve_{J'})+\rho_I(\ve_J-\ve_{J'})\right)=\dfrac12(\ve_J+\ve_{J'}-(\ve_J-\ve_{J'}))=\ve_{J'}.$$

For the second part, suppose $K\in\mathcal{C}(J).$ Then $K'=w_J(K)\subset J\subset I$ and  
$$\rho_I(\ve_K-\ve_{K'})=\rho_I(\ve_K)-\rho_I(\ve_{K'})=\ve_{w_I(K)}-\ve_{w_I(w_J(K))}.$$
Note that $w_Iw_J=w_Iw_Jw_I^{-1}w_I=w_{J'}w_I.$ Then
$$\rho_I(\ve_K-\ve_{K'})=\ve_{w_I(K)}-\ve_{w_I(w_J(K))}=\ve_{w_I(K)}-\ve_{w_{J'}(w_I(K))}\in E_{J'}$$
since $w_I(K)\subset w_I(J)=J'.$ Hence $\rho_I(E_J)\subset E_{J'}.$ A similar argument shows that $\rho_I(E_{J'})\subset E_J$ and then  $\rho_I(E_J)=E_{J'}$ as $\rho_I$ is an involution. 
\end{proof}
Case (b) now follows from the commutative diagram implied by part (1) of the lemma:
\begin{center}
\begin{tikzpicture}
	\matrix (m)
		[
			matrix of math nodes,
			row sep    = 3em,
			column sep = 4em
		]
		{
			(\mathbb R\ve_J\oplus E_J)\oplus F_J   & (\mathbb R\ve_J\oplus E_J)\oplus F_J \\
			(\mathbb R\ve_J\oplus E_J)\oplus F_J   & (\mathbb R\ve_J\oplus E_J)\oplus F_J \\
		};
	\path
		(m-1-1) edge [->] node [left] {$\rho_J=-\mathrm{id}\oplus \mathrm{id}$} (m-2-1)
		(m-1-1) edge [->] node [above] {$\rho_I$} (m-1-2)
		(m-2-1) edge [->] node [above] {$\rho_I$} (m-2-2)
		(m-1-2) edge [->] node [right] {$\rho_J=-\mathrm{id}\oplus \mathrm{id}$} (m-2-2);
\end{tikzpicture}
\end{center}

Part (2) of the lemma gives a similar commutative diagram for case (c):
\begin{center}
\begin{tikzpicture}
	\matrix (m)
		[
			matrix of math nodes,
			row sep    = 3em,
			column sep = 4em
		]
		{
			(\mathbb R\ve_J\oplus E_J)\oplus F_J   & (\mathbb R\ve_{J'}\oplus E_{J'})\oplus F_{J'} \\
			(\mathbb R\ve_J\oplus E_J)\oplus F_J   & (\mathbb R\ve_{J'}\oplus E_{J'})\oplus F_{J'} \\
		};
	\path
		(m-1-1) edge [->] node [left] {$\rho_J=-\mathrm{id}\oplus \mathrm{id}$} (m-2-1)
		(m-1-1) edge [->] node [above] {$\rho_I$} (m-1-2)
		(m-2-1) edge [->] node [above] {$\rho_I$} (m-2-2)
		(m-1-2) edge [->] node [right] {$\rho_{J'}=-\mathrm{id}\oplus \mathrm{id}$} (m-2-2);
\end{tikzpicture}
\end{center}

This completes the proof.
\end{proof}
When $t=0,$ $\ve_J$ is orthogonal to $\ve_I$ and $E_I$ and $\rho(\ve_J)=\ve_J.$ Lemma 5 implies the images of the standard basis $\{\ve_J\}_{J\in\mathcal F(S)}$ under any $\rho(w)$ for $w\in C_W$ lie in the set $\{\pm\ve_J\}_{J\in\mathcal F(S)}$ when $t=0.$ Hence the representation $\rho:C_W\to GL(E)$ cannot be faithful when $C_W$ is infinite, nor will $\rho$ be faithful when restricted to an infinite subgroup of $C_W$ (for example the generalized pure cactus group $PC_W$).

\subsection{Examples of Generalized Cactus Groups}
In this section, we analyze a few examples of generalized cactus groups and their representations.
\subsubsection{Type $A_2$} 
The cactus group $J_3$ is defined on the Coxeter group $S_3$ and is generated by $s_1$ and $s_2$ with $$\mathcal F=\{\{s_1\},\{s_2\},\{s_1,s_2\}\}.$$ Then $J_3$ can be generated by $B=\gamma_{\{s_2\}}$ and $C=\gamma_{\{s_1,s_2\}}$ and the representation $\rho:J_3\to GL(\mathbb R,3)$ is given by involutions
$$
\rho(B)=\begin{pmatrix}
	1 &  0 & 0 \\
	2t & -1 & 0 \\
	0 & 0 & 1
\end{pmatrix},
\rho(C)=\begin{pmatrix}
	0 & 1 & 0 \\
	1 & 0 & 0 \\
	0 & 0 & -1
\end{pmatrix}
$$
and then $\left.\rho\right|_U$ is irreducible where $U=\mathrm{span}(\ve_1,\ve_2).$ 

Furthermore, the group $J_3$ embeds into $D_3\rtimes \aut D_3=D_3\rtimes S_3$ by Theorem \ref{embed}. Then the representation $\Pi:J_3\hookrightarrow D_3\rtimes S_3\to GL(\mathbb R,4)$ is given by
$$
\Pi(B)=\begin{pmatrix}
	0 & 1 & 0 & 0 \\
	1 & 0 & 0 & 0 \\
	2t & 2t & -1 & 0 \\
	0 & 0 & 0 & 1
\end{pmatrix},\,
\Pi(C)
=\begin{pmatrix}
	0 & 0 & 1 & 0 \\
	0 & 1 & 0 & 0 \\
	1 & 0 & 0 & 0 \\
	0 & 0 & 0 & -1
\end{pmatrix}.
$$
The subspace $U_1=\mathrm{span}(\ve_1,\ve_2,\ve_3)$ is invariant under $\Pi$ and a stable vector of $\Pi$ is $v_0=\ve_1-\ve_2+\ve_3.$ Hence there is a quotient representation of $\Pi$ on the quotient space $U_1/\langle v_0\rangle$ with basis $\bar \ve_1$ and $\bar \ve_3.$ Under this basis, the quotient representation is given by 
$$
B\mapsto\begin{pmatrix}
	1 & 0 \\
	2t+1 &  -1
\end{pmatrix},\,
C\mapsto\begin{pmatrix}
	0 & 1 \\
	1 & 0
\end{pmatrix}.
$$
The quotient representation is isomorphic to the representation $\rho$ with a change in parameter $t$ to $t+1/2.$
\subsubsection{Type $I_2(n)$}
The dihedral group of order $n$ is generated by $\{a,b\}$ with an associated Coxeter matrix 
$$
\begin{pmatrix}
	1 & n \\
	n & 1
\end{pmatrix}
$$
and $\mathcal F=\{\{a\},\{b\},\{a,b\}\}.$ The longest element in $D_{2n}$ is $(ab)^{n/2}$ if $n$ is even and $b(ab)^{(n-1)/2}$ if $n$ is odd.

Let $CD_{2n}$ be the generalized cactus group on $D_{2n}$ generated by $\gamma_1=\gamma_a,\gamma_2=\gamma_b$ and $\gamma_3=\gamma_{ab}.$ Then $CD_{2n}$ is isomorphic to $(\mathbb Z_2)^3$ when $n=2,$ to $D_{\infty}\times \mathbb Z_2$ when $n\ge 4$ and $n$ is even, and to $J_3$ (hence to $D_\infty$) when $n$ is odd.

Let $n$ be an odd integer, and then $CD_{2n}$ admits the same representation $\rho_t$ on $E$ as $J_3$ does. However, the basis of the vector space $\mathbb E$ consists of the vectors $\ve_i$ for $0\le i\le n,$ indexed by the subgroups $A_i=\{e,(ab)^ia\}$ for $i=0,\cdots,n-1$ (as $b=(ab)^{n-1}a$ is a conjugation of $a$) and $D_{2n}$ for $i=n$. Using the notation in the previous section where $B=\gamma_b$ and $C=\gamma_{ab},$ we have
$$\Pi(B)(\ve_i)=\begin{cases}
	\ve_{n-i-2}+2t\ve_{n-1},&0\le i\le n-2 \\
	-\ve_{n-1}, &i=n-1 \\
	\ve_n,&i=n
\end{cases}$$
and
$$
\Pi(C)(\ve_i)=\begin{cases}
	\ve_{n-1-i}, & 0\le i\le n-1 \\
	-\ve_n, &i=n.
\end{cases}
$$

There are two subrepresentations of $\Pi$ spanned by $\ve_i$ for $0\le i\le n-1$ and $\ve_n,$ respectively. A stable subspace $V$ of $\Pi$ is spanned by $\ve_i-\ve_{n-2-i}$ for $0\le i\le n-2,$ and a stable vector of $\Pi$ is 
$$v=\ve_0-\ve_1+\cdots+\ve_{2k}.$$

\section*{Acknowledgements}
The author would like to warmly thank Professor Raphaël Rouquier for supporting this work. This project was partly funded by NSF Grant DMS-1702305 as part of the 2021 Summer REU Program at UCLA.

\nocite{*}
\printbibliography

\end{document}